\documentclass{amsart}
\usepackage[dvips]{graphicx}

\newtheorem{theorem}{Theorem}[section]
\newtheorem{lemma}[theorem]{Lemma}
\newtheorem{proposition}[theorem]{Proposition}
\newtheorem{corollary}[theorem]{Corollary}

\theoremstyle{definition}
\newtheorem{definition}[theorem]{Definition}
\newtheorem{remark}[theorem]{Remark}

\numberwithin{equation}{section}

\newcommand{\zz}{\mathbb{Z}}

\newcommand{\cpk}{\mathbb{CP}^2}
\newcommand{\cpkk}{\overline{\mathbb{CP}}{}^2}

\begin{document}

\title[Homologous Non-isotopic Symplectic Tori
in Homotopy $E(1)$'s]{Homologous Non-isotopic
Symplectic Tori\\
in Homotopy Rational Elliptic Surfaces}

\author{Tolga Etg\"u}
\address{Department of Mathematics and Statistics, McMaster University,
Hamilton, Ontario L8S 4K1, Canada}
\email{etgut@math.mcmaster.ca}

\author{B. Doug Park}
\address{Department of Pure Mathematics, University of Waterloo, Waterloo,
Ontario, N2L 3G1, Canada}
\email{bdpark@math.uwaterloo.ca}
\thanks{B.D. Park was partially supported by an NSERC research grant.}

\date{May 15, 2003.  Revised on \today}

%----------------------------------------------------------------------------

\begin{abstract}
Let $E(1)_K$ denote the homotopy rational elliptic surface corresponding to a
knot $K$ in $S^3$ constructed by
R. Fintushel and R.J. Stern in \cite{fs:knots}.  We
construct an infinite family of homologous non-isotopic symplectic
tori representing a primitive $2$-dimensional
homology class in $E(1)_K$ when
$K$\/ is any nontrivial fibred knot in $S^3$.
We also show how these tori can be non-isotopically embedded as
homologous symplectic submanifolds in other symplectic $4$-manifolds.
\end{abstract}

\maketitle

%------------------------------------------------------------------------------

\section{Introduction}

This paper is a continuation of studies initiated in \cite{ep1}
and \cite{ep:k3} regarding infinite families of non-isotopic and symplectic
tori representing the same homology class in a symplectic 4-manifold.
Let $E(1)_K$ denote the closed 4-manifold that is homotopy equivalent
(hence homeomorphic) to the rational elliptic surface $E(1)\cong\cpk \#9\cpkk$
and is obtained by performing knot surgery (as defined in \cite{fs:knots})
on the rational elliptic surface
using a knot $K$\/ in $S^3$. Our main result is the following:

\begin{theorem}\label{theorem:main}
Let $K\subset S^3$ be a nontrivial fibred knot. Then there
exists an infinite family of pairwise non-isotopic symplectic tori
representing the primitive homology class\/ $[F]=[T_m]$ in $E(1)_K$,
where\/ $[F]$ is the homology
class of the fiber in a rational elliptic surface $E(1)\cong \cpk \#
9\cpkk$.
\end{theorem}

Examples of homologous, non-isotopic, symplectic tori were first
constructed in \cite{fs:non-isotopic} and then in \cite{ep1},
\cite{ep:k3}, \cite{vidussi:non-isotopic} and
\cite{vidussi:E(1)_K} (also see \cite{fs:lagrangian} and
\cite{vidussi:lagrangian}).  Recall that infinite families of
non-isotopic symplectic tori representing $n[F]\in H_2(E(1)_K)$,
$n\geq 2$, were constructed in \cite{ep1}.  The family of tori we
construct in this paper is in some sense the `simplest' example
known so far, when measured in terms of the `geography size' of
the ambient (simply-connected) symplectic 4-manifold, the
divisibility of the homology class represented, and the complexity
of the knotting of the tori.  In \cite{vidussi:E(1)_K}, using a
different technique, Vidussi already constructed symplectic tori
representing the same primitive class in $E(1)_K$ for some
particular fibred $K$, namely the trefoil and other fibred knots
that have the trefoil as one of their connected summands.

It should be noted that the
non-existence of such an infinite family of tori in $\cpk$ and
$\cpk \# \cpkk$ is proved by Sikorav in \cite{sikorav} and by Siebert and Tian
in \cite{st}, respectively. It is also conjectured that there is at most
one
symplectic torus (up to isotopy) representing each homology class in
$\cpk \# n \cpkk$ for
$n < 9$.

The proof of Theorem~\ref{theorem:main} is spread out over the next three
sections. We will review the relevant definitions in
Section~\ref{sec:link surgery}. In Section~\ref{sec:generalization},
we will present a
direct generalization in the form of Proposition~\ref{prop:generalization}.
In this introduction and elsewhere in the paper by isotopy we mean smooth
isotopy
and all homology groups have $\zz$ coefficients.

%--------------------------------------------------------------------------

\section{Link Surgery 4-Manifolds}
\label{sec:link surgery}

In this section, first we review the generalization of the
link surgery construction of
Fintushel and Stern \cite{fs:knots} by Vidussi
\cite{vidussi:smooth}, and then give specific link surgeries that
will be used in the following sections.

For an $n$-component link $L\subset S^3$,
choose an ordered homology basis of simple closed curves $\{(\alpha_i, \beta_i)
\}_{i=1}^{n}$ such that
the pair $(\alpha_i, \beta_i )$ lie in the $i$-th boundary
component of the link exterior
and the intersection of
$\alpha_i$ and $\beta_i$ is 1.
Let $X_i$
($i=1,\dots, n$) be a 4-manifold containing a 2-dimensional torus
submanifold $F_i$ of self-intersection $0$.  Choose a Cartesian product
decomposition $F_i = C_1^{i} \times C_2^{i}$, where each $C^i_j \cong S^1$
($j=1,2$) is an embedded
circle in $X_i$.

\begin{definition}\label{def:data}
The ordered collection
$$\mathfrak{D} \, =\; \big( \{(\alpha_i, \beta_i)
\}_{i=1}^{n}\: , \: \{(X_i  , \hspace{1pt}
F_i=  C_1^{i} \times C_2^{i} )\}_{i=1}^{n}\big)$$
is called a \emph{link
surgery gluing data}\/ for an $n$-component link $L$.
We define the \emph{link surgery manifold corresponding to} $\mathfrak{D}$
to be the closed $4$-manifold
\[
L(\mathfrak{D}) \: :=\; [\coprod_{i=1}^{n} X_i\setminus\nu
F_i]\hspace{-20pt}\bigcup_{F_i\times\partial D^2=(S^1\times
\alpha_i)\times\beta_i}\hspace{-20pt} [S^1\times(S^3\setminus \nu
L)]\, ,
\]
where $\nu$\/ denotes the tubular neighbourhoods.  Here, the gluing
diffeomorphisms between the boundary 3-tori identify the torus
$F_i =  C_1^{i} \times C_2^{i}$\/ of $X_i$\/
with\/ $S^1\times \alpha_i$\/ factorwise, and act as the complex
conjugation on the last remaining\/ $S^1$ factor.
\end{definition}

\begin{remark}
Strictly speaking, the diffeomorphism type of the link surgery manifold
$L(\mathfrak{D})$\/ may possibly
depend on the chosen trivialization of $\,\nu F_i \cong F_i \times D^2$ (the
framing of $F_i$).  However, we will suppress this dependence in our notation.
It is well known (see e.g. \cite{gs}) that the diffeomorphism type of
$L(\mathfrak{D})$
is independent of the framing of $F_i$ when $(X_i,F_i) = (E(1), F)$.
\end{remark}

We fix a Cartesian product decomposition\/
$F= C_1 \times C_2$ in $E(1)$.
Let $K$\/ be a knot in $S^3$, and
let $M_K$ denote the 3-manifold that is the result
of the 0-framed surgery on $K$.
Fix a meridian circle $m=\mu(K)$\/ in $M_K$.

\begin{figure}[!ht]
\begin{center}
\includegraphics[scale=.5]{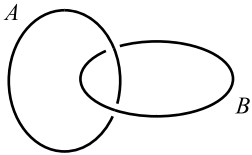}
\end{center}
\caption{Hopf link $L= A\cup B$}
\label{fig:hopf}
\end{figure}

\begin{definition}
Let\/ $L \subset S^3$ be the Hopf link in Figure~$\ref{fig:hopf}$.
For the link surgery gluing data\/
\begin{eqnarray}\label{eq:data}
 \mathfrak{D} \!\!\!  &:=& \!\!\!
\big(\{ (\mu(A),\lambda(A)), (\lambda(B),-\mu(B))
\}, \\ && \{ (X_1, F_1=C_1^1\times C_2^1),
(S^1\times M_K, T_m = S^1\times m )\}\big),
\nonumber
\end{eqnarray}
we shall denote\/ $L(\mathfrak{D})$ by $(X_1)_K$.
Here, $\mu(\,\cdot\,)$ and $\lambda(\,\cdot\,)$ denote
the meridian and the longitude
of a knot, respectively.
In particular, when $(X_1, F_1=C_1^1\times C_2^1) = (E(1), F=C_1\times C_2)$,
we denote $L(\mathfrak{D})$\/ by $E(1)_K$. This notation is consistent
with that of Fintushel and Stern in \cite{fs:knots} as
there is a diffeomorphism
between our $L(\mathfrak{D})$ and their fiber sum
$E(1)_K = E(1) \#_{F=T_m}\! (S^1 \times M_K)$.
\end{definition}

Note that there is a canonical framing of $T_m$ in $(S^1\times M_K)$\/ given by
the minimal genus
Seifert surface of the knot $K$.  We shall always use this framing to
trivialize $\nu T_m$.

\begin{lemma}\label{lemma:E(1)_K is symplectic}
If\/ $K\subset S^3$ is a fibred knot, then $E(1)_K$ is a symplectic
$4$-manifold.
\end{lemma}

\begin{proof}
This is because there exists a fiber bundle
$\,(S^1\times M_K) \rightarrow T^2$ when $K$\/ is
fibred, so
$(S^1\times M_K)$ admits a symplectic form with respect to which
$T_m$ is a symplectic submanifold (cf.$\;$\cite{thurston}).
Hence we may express $E(1)_K$ as a \emph{symplectic}\/ fiber sum
$E(1)\#_{F=T_m}\! (S^1\times M_K)$ along symplectic submanifolds $F$ and $T_m$
(cf.$\;$\cite{g:sum}).
\end{proof}

\begin{lemma}
The homology class\/ $[F]=[T_m]\in H_2(E(1)_K)$
is primitive.
\end{lemma}

\begin{proof}
Since\/ $[\mu(A)]=[\lambda(B)] \in H_1(S^3\setminus \nu L)$,
we must have\/ $[S^1\times \mu(A)]=[S^1\times \lambda(B)]$\/ in\/
$H_2(S^1\times (S^3\setminus \nu L) )$, and so\/
$[F]=[T_m]$\/ in\/ $H_2(E(1)_K)$.
Let $\Sigma$ denote a Seifert surface of $K$.
Let $\Sigma_K$ denote a closed surface in $E(1)_K$
that is the internal tubular sum of
a punctured section in $[E(1)\setminus \nu F]$ and a punctured surface
$\{ {\rm point} \} \times\Sigma$, glued together along $K$. Then we have
$[\Sigma_K] \cdot [T_m]= \pm 1$.
\end{proof}

%---------------------------------------------------------------------------

\section{Family of Homologous Symplectic Tori in $E(1)_K$}
\label{sec:symplectic family}

\begin{figure}[!ht]
\begin{center}
\includegraphics[scale=.85]{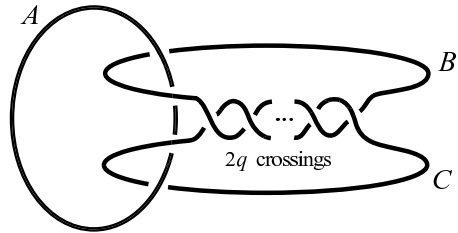}
\end{center}
\caption{3-component link $L_q=A\cup B \cup C\,$ in $S^3$}
\label{fig:doublehopf}
\end{figure}

Let $T_C := S^1 \times C \subset [S^1 \times (S^3\setminus \nu L)]
\subset E(1)_K$,
where the closed curve $C\subset (S^3\setminus \nu L )$\/
is given by Figure~\ref{fig:doublehopf}.

\begin{lemma}
If\/ $K$ is a fibred knot, then $T_C =S^1\times C$ is a
symplectic submanifold of\/ $E(1)_K$ and we have\/
$[T_C]=[F]\,$ in\/ $H_2(E(1)_K)$.
\end{lemma}

\begin{proof}
It is easy to see that the link exterior
$Y:=(S^3\setminus \nu L)$\/ is diffeomorphic to
$S^1 \times \mathbb{A}$, where $\mathbb{A}\cong S^1 \times [0,1]\,$
is an annulus.
Hence we have
\[
[S^1 \times (S^3\setminus \nu L)] \:\cong\:
[S^1 \times (S^1 \times \mathbb{A})]
\:\cong\: T^3\times [0,1] .
\]
We may assume that the symplectic form on $E(1)_K$ restricts to
\[
\omega \,=\, dx\wedge dy \,+\, r\,dr\wedge d\theta
\]
on $[S^1 \times (S^1\times \mathbb{A})]$,
where $x$\/ and $y$\/ are the angular coordinates on the first and the second
$S^1$ factors respectively, and
$(r,\theta)$ are the polar coordinates on the annulus $\mathbb{A}$.
We can embed the curve $C$\/ inside $(S^1\times \mathbb{A})$ such that
$C$ is transverse to every annulus of the form, $\{{\rm point}\} \times
\mathbb{A}$, and the restriction $dy|_{C}$ never vanishes.
It follows that $\omega |_{T_C} = (dx\wedge dy)|_{T_C} \neq 0$, and
consequently
$T_C$ is a symplectic submanifold of $E(1)_K$.

To determine the homology class of $T_C$, note that\/ $[C]=[\mu
(A)] + q [\mu (B)]$\/ in\/ $H_1(Y)$.
When we glue\/ $[(S^1\times M_K)\setminus \nu T_m]$\/ to\/  $S^1\times
Y$, the homology class $[\mu(B)]$ gets identified with
$[\{{\rm point}\}\times\lambda(K)]
\in H_1((S^1\times M_K)\setminus \nu T_m)$, which is
trivial.  Hence by K\"unneth's theorem, $[T_C]=[S^1\times
\mu(A)]\,$ in $H_2([S^1\times Y]\cup [(S^1\times M_K)\setminus \nu
T_m])$. It follows that\/ $[T_C]=[F]\,$ in $H_2(E(1)_K)$.
\end{proof}

%------------------------------------------------------------------------------

\section{Non-Isotopy:  Seiberg-Witten Invariants}
\label{sec:sw}

Our strategy is to show that the isotopy types of the tori $\{ T_C
\}_{q\geq 1}$ can be distinguished by comparing the Seiberg-Witten
invariants of the corresponding family of fiber sum 4-manifolds $\{E(1)_K
\#_{T_C=F} E(n)\}_{q\geq 1, n\geq 1}$. Note that there is a
canonical framing of a regular fiber $F$\/ in $E(n)$, coming from
the elliptic fibration $E(n)\rightarrow \mathbb{CP}^1$.

\begin{lemma}
The fiber sum\/ $E(1)_K \#_{T_C=F} E(n)$\/ is diffeomorphic to the
link surgery manifold $L_q(\mathfrak{D}')$, where
\begin{eqnarray}
&&  \quad\quad   \mathfrak{D}' : = \:  \big(
\{(\mu(A),\lambda(A)),(\lambda(B),-\mu(B)),
(\lambda(C),-\mu(C))\},  \\[3pt]
&& \{(E(1),F=C_1\times C_2),
(S^1\times M_K,T_m=S^1 \times\mu(K)),(E(n),F=C_1\times C_2)\}\big).
\nonumber
\end{eqnarray}
\end{lemma}

\begin{proof}
We already observed in the proof of
Lemma~\ref{lemma:E(1)_K is symplectic} that the fiber sum construction
corresponds to this type of link surgery.  (See also \cite{ep:k3}.)
\end{proof}

Recall that the Seiberg-Witten invariant\/
$\overline{SW}_{\!\!X}$\/ of a 4-manifold $X$\/
%(satisfying\/ $b_2^+(X)>1$)
can be thought of as an element of the group ring of $H_2(X)$,
i.e.\/ $\overline{SW}_{\!\!X} \in \zz [ H_2( X ) ]$. If we
write\/ $\overline{SW}_{\!\!X} = \sum_g a_g g \hspace{1pt}$, then
we say that\/ $g\in H_2( X )$\/ is a Seiberg-Witten
\emph{basic class}\/ of $X$\/ if\/ $a_g\neq 0$. Since the
Seiberg-Witten invariant of a 4-manifold is a diffeomorphism
invariant, so are the divisibilities of Seiberg-Witten basic
classes. The Seiberg-Witten invariant of the link surgery manifold
$L_q(\mathfrak{D}')$ is known to be related to the Alexander
polynomial\/ $\Delta_{L_q}$ of the link $L_q$.

\begin{lemma}
$\Delta_{L_q}(x,s,t)= 1-x(st)^q$, where the variables\/ $x$, $s$
and\/ $t$\/ correspond to the components\/ $A$, $B$\/ and\/ $C$\/
respectively.
\end{lemma}

\begin{proof}
This follows readily from the formula in Theorem 1 of \cite{morton}
which gives the multivariable Alexander polynomial
of a braid closure and its axis
in terms of the representation of the braid.
We view $A$\/ as the axis of the closure of a 2-strand braid, $B\cup C$.
See \cite{ep1} for details on a similar computation.
\end{proof}

\begin{theorem}\label{theorem:seiberg-witten}
Let\/ $\iota:[S^1\times(S^3\setminus\nu L_{q})]\rightarrow
L_{q}(\mathfrak{D}')$\/ be the inclusion map.  Let\/
$\xi:=\iota_{\ast}[S^1\times\mu(A)],$
$\tau:=\iota_{\ast}[S^1\times\mu(C)] \in
H_2(L_{q}(\mathfrak{D}') ).$  Then\/ $\xi$ and\/ $\tau$ are both primitive
and linearly independent.
The Seiberg-Witten invariant of\/ $L_{q}(\mathfrak{D}')$ is
given by
\begin{equation}\label{eq:sw}
\overline{SW}_{\!L_q(\mathfrak{D}')}\, =  \: (\xi^{-1}-\xi)^{n-1}
\cdot \Delta^{\rm sym}_K(\xi^2 \tau^{2q})\hspace{1pt} ,
\end{equation}
where\/ $\Delta^{\rm sym}_K$ is the\/ \emph{symmetrized}
Alexander polynomial of the knot\/ $K$.
\end{theorem}

\begin{proof}
Let\/ $N := (S^3\setminus \nu L_q)$, and let $Z:=
[(S^1\times M_K) \setminus \nu T_m]$. Recall from \cite{doug:pft3}
that we have $\overline{SW}_{\!E(n)\setminus \nu
F}=([F]^{-1}-[F])^{n-1}$, and also
$$
\overline{SW}^{\,\pm}_{\!Z} = \:
\overline{SW}^{\,\pm}_{\!(S^1\times M_K) \setminus \nu T_m} =\:
\frac{\Delta_K^{\rm sym}([T_m]^2)}{[T_m]^{-1}-[T_m]}
\; .$$
From the gluing formulas in \cite{doug:pft3} and \cite{Taubes:T^3},
we may conclude that
$$
\overline{SW}_{\!L_q(\mathfrak{D}')} \,=\:
\overline{SW}_{\!E(1)\setminus \nu F} \cdot
\overline{SW}_{\!E(n)\setminus \nu F} \cdot
\overline{SW}^{\,\pm}_{\!(S^1\times M_K) \setminus \nu T_m} \cdot
\Delta_{L_q}^{\rm sym}(\xi^2,\sigma^2,\tau^2),
$$
where $\sigma:=\iota_{\ast}[S^1\times\mu(B)]$. Note that $\sigma
=1\in \zz[H_2(L_q(\mathfrak{D}'))]$,
since we have\/ $[\mu(B)]=[\lambda(K)]=0\,$ in\/ $H_1(Z)$.
Also note that $\mu(K)$ and $\lambda(B)$ are identified
by the gluing data $\mathfrak{D}'$, and\/
$[\lambda(B)]=[\mu(A)]+q[\mu(C)]\in H_1(N)$.  It follows
that\/ $[T_m]=\iota_{\ast}[S^1\times \lambda(B)]=
\xi \tau^q \in \zz[H_2(L_q(\mathfrak{D}'))]$. Thus
we have
$$\Delta_{L_q}^{\rm sym}(\xi^2,\sigma^2,\tau^2) \:=\:
\xi^{-1}\tau^{-q} - \xi \tau^q \:=\: [T_m]^{-1}-[T_m]\, .$$
Hence
\begin{equation}\label{eq:sw almost}
\overline{SW}_{\!L_q(\mathfrak{D}')}
\,=\:([F]^{-1}-[F])^{n-1}\cdot \Delta_K^{\rm sym}(\xi^2\tau^{2q})\hspace{1pt} .
\end{equation}
Note that the fiber $F$\/ in\/ $E(n)$\/ gets identified with\/
$S^1\times \lambda(C)$ by the gluing data $\mathfrak{D}'$,
and we also have\/ $[\lambda(C)]=[\mu(A)]
+q[\mu(B)]=[\mu(A)]\,$ in\/ $H_1([S^1\times N]\cup Z )$.
Therefore we can identify\/ $[F]=\xi\,$ in (\ref{eq:sw almost}),
and we obtain Equation
(\ref{eq:sw}).

Next we show that $\xi$ and $\tau$ are primitive and linearly
independent elements of $H_2(L_{q}(\mathfrak{D}'))$. We can
proceed in two different ways.  A Mayer-Vietoris argument,
combined with Freedman's classification theorem
(cf.$\;$\cite{FQ}), shows that $L_q(\mathfrak{D}')$\/ is
homeomorphic to $E(n+1)$. It is not too hard to find two closed
surfaces $R$\/ and $S$\/ in $L_q(\mathfrak{D}')$\/ satisfying
\[
\xi \cdot [S] \;=\; \tau \cdot [R] \;=\;  1 \, ,
\]
and
\[
\xi \cdot [R] \;=\; \tau \cdot [S] \;=\; [R] \cdot [S] \;=\; 0 \, .
\]
For example, we can let $S$\/ be the
internal tubular sum of punctured sections from\/
$[E(1)\setminus \nu F]$\/ and\/ $[E(n)\setminus \nu F]$\/ summands,
together with a suitable punctured surface from the $Z$\/ summand.
Let $R$\/ be
the internal tubular sum of the self-intersection $(-1)$ disks bounding
the circle $C_2$ in
$[E(1)\setminus \nu F]$\/ and\/ $[E(n)\setminus \nu F]$\/ summands,
together with a suitable punctured surface from the $Z$\/ summand.

In $L_q(\mathfrak{D}')$, $S$\/ plays the role of a section in
$E(n+1)$, while $\xi$\/ plays the role of the homology class of
the fiber. Note that we have\/ $[\mu(C)]=[\lambda (A)]-[\mu(B)]
=[\lambda(A)]\in H_1([S^1\times N]\cup Z)$, and
the gluing data $\mathfrak{D}'$ identifies
$\lambda(A)$
with a meridian circle $\mu(F)$ of the fiber $F$\/ in
$\partial [E(1)\setminus \nu F]$. Hence $\tau$ plays the role of the
homology class of the rim torus\/ $C_1 \times \mu (F)$\/ in $E(n+1)$.
$R$\/ plays the role of a self-intersection $(-2)$ sphere transversally
intersecting the above rim torus once.  The pairs\/ $(\xi ,[S])$\/
and\/ $(\tau, [R])$\/ form homology bases for two $\bigl(
\begin{smallmatrix}
  0 & 1 \\
  1 & 0
\end{smallmatrix} \bigr)$ summands in the intersection form of
$L_q(\mathfrak{D}')$.

Alternatively, we can argue more algebraically as follows.
Consider the composition of homomorphisms
\begin{equation}\label{eq:homomorphism}
H_1(N) \longrightarrow H_2(S^1\times N)
\stackrel{\iota_{\ast}\;}{\longrightarrow}
H_2(L_q(\mathfrak{D}')),
\end{equation}
where the first
map is a part of the K\"unneth isomorphism
\begin{equation}\label{eq:Kunneth}
H_1(N) \oplus H_2(N)
\stackrel{\cong}{\longrightarrow} H_2(S^1\times N).
\end{equation}
Note that $H_2(N) \cong \zz\oplus \zz\,$,
as is easily seen from the long exact
sequence of the pair\/ $(N,\partial N)$ as follows.
\[
\begin{array}{ccccccccc}
H^0(N)&\longrightarrow & H^0(\partial N) &\longrightarrow &
 H^1(N,\partial N) & \stackrel{0}{\longrightarrow} &
H^1(N) & \longrightarrow &  H^1(\partial N) \\
_{||} && _{||} && _{||} & & _{||}
&& _{||}\\[3pt]
\zz  & \longrightarrow  & \zz^3 & \longrightarrow & H_2(N) &
 \stackrel{0}{\longrightarrow}
& \zz^3 &
\longrightarrow & \zz^6
\end{array}
\]
Note that the first map sends the generator $\hspace{1pt} 1\in \zz\,$ to
the diagonal element\/ $(1,1,1)\in \zz^3$, while the last map is injective.
We have also used the Lefschetz duality theorem (for manifolds with boundary)
to identify $H_2(N)\cong
H^1(N,\partial N)$.

Next consider the long exact sequence of the pair
$(L_q(\mathfrak{D}'), S^1\times N)$:
\[
0 = H_3(L_q(\mathfrak{D}')) \longrightarrow
H_3(L_q(\mathfrak{D}'),S^1\times N) \longrightarrow H_2(S^1\times
N) \stackrel{\iota_{\ast}\;}{\longrightarrow}
H_2(L_q(\mathfrak{D}'))
\]
The kernel of the last map $\iota_{\ast}$ is isomorphic to
$H_3(L_q(\mathfrak{D}'),S^1\times N)$.  By Lefschetz duality
theorem (for relative manifolds),
$H_3(L_q(\mathfrak{D}'),S^1\times N)$ is in turn
isomorphic to
\[
H^1(L_q(\mathfrak{D}')\setminus (S^1\times N)) \:\cong\;
H^1(E(1)\setminus \nu F) \oplus H^1(E(n)\setminus \nu F) \oplus
H^1(Z).
\]
Since we have $H^1(E(1)\setminus \nu F) = H^1(E(n)\setminus \nu F)=0\,$
and
\begin{equation}\label{eq:z plus z}
H^1(Z)
=H^1(S^1 \times( M_K \setminus \nu m) ) \cong \; \zz\oplus \zz \, ,
\end{equation}
the kernel of $\iota_{\ast}$ is isomorphic to\/ $\zz\oplus\zz\,$.

Finally we observe that only one $\zz$ summand of (\ref{eq:z plus z}) lies
in the kernel of the composition (\ref{eq:homomorphism}).  The other
$\zz$ summand belongs to the kernel of
\[
H_2(N) \longrightarrow H_2(S^1\times N)
\stackrel{\iota_{\ast}\;}{\longrightarrow}
H_2(L_q(\mathfrak{D}')),
\]
where the first map is the second part of the K\"unneth
isomorphism (\ref{eq:Kunneth}). We have thus shown that the kernel
of the composition (\ref{eq:homomorphism}) is of rank one.  It
follows immediately that $\xi$ and $\tau$ are linearly
independent, since we already have shown that $\sigma$ is trivial.
A more detailed analysis shows that\/ $\{\xi ,\tau\}$\/ can be
extended to a basis of $H_2(L_q(\mathfrak{D}'))$, which we shall
omit. (Also see the proof of Proposition 3.2 in
\cite{McMullen-Taubes} for a similar argument.)
\end{proof}

\begin{corollary}\label{cor:divisibilities}
If\/ $K$ is a nontrivial fibred knot,
then the tori\/ $\{ T_C \}_{q \geq 1}$ are
pairwise non-isotopic inside $E(1)_K$.  In fact, there is no
self-diffeomorphism of\/ $E(1)_K$ that maps one element of this
family to another.
\end{corollary}

\begin{proof}
Let's choose $n$\/ to be $2g+1$, where $g$ is the genus of $K$.
Remember that the degree of the symmetrized Alexander polynomial of a fibred
knot is the same as its genus (see e.g. Proposition 8.16 in \cite{bz}).
 Since we assume that $K$ is nontrivial, i.e. not the unknot,
$g >0$. A Seiberg-Witten basic class of $L_q(\mathfrak{D}')$ with the
highest divisibility is divisible by $2gq$.
This could be seen by observing that
the highest power of $\tau$ in (\ref{eq:sw}) of
Theorem~\ref{theorem:seiberg-witten} is $2gq$
(hence there cannot be a basic class
with divisibility higher than $2gq$) and our choice of $n=2g+1$\/ ensures
that there is a basic class (namely $\tau^{2gq}$) with this highest possible
divisibility. On the other hand, since the
Seiberg-Witten invariant is a diffeomorphism invariant, so are
the divisibilities of basic classes. Therefore, $L_q(\mathfrak{D}')$
is diffeomorphic to $L_{q'}(\mathfrak{D}')$ if and only if\/ $q=q'$
proving that the tori
in $\{ T_C \}_{q \geq 1}$ are different up to isotopy and in fact even up to
self-diffeomorphisms of\/ $E(1)_K$.
\end{proof}

This concludes the proof of Theorem~\ref{theorem:main}.

%----------------------------------------------------------------------------

\section{Generalization to Other Symplectic 4-Manifolds}
\label{sec:generalization}

When $K$\/ is the unknot, $E(1)_K$ is diffeomorphic to $E(1)$. In
this unknot case, our family of tori $\{T_C\}_{q \geq 1}$ are
easily seen to be all isotopic to one another.  The isotopy can
actually be visualized by erasing the $B$\/ component in
Figure~\ref{fig:doublehopf} (This corresponds to filling in $\nu
B$\/ with\/ $(M_K\setminus \nu m)$, which, in the unknot case, is
diffeomorphic to a solid torus\/ $S^1\times D^2$.), and
straightening out the $C$\/ component through the tubular
neighbourhood of $B$, which has now been filled in. Note that the
normal disks of $B$\/ are the Seifert surfaces of the unknot.

Suppose that $K$\/ is not fibred. Then, unlike the fibred case where the degree
of the Alexander polynomial is (the same as the genus hence)
strictly greater than
$0$ unless the knot is the unknot, the Alexander polynomial of $K$\/ might
be constant and the Seiberg-Witten invariant doesn't seem to be
delicate enough to distinguish the tori we constructed. On the
other hand, when $K$ is not fibred and has a non-constant Alexander
polynomial, the tori in our family\/ $\{ T_C \}_{q\geq 1}$ are still pairwise
non-isotopic in $E(1)_K$, but there is no natural symplectic
structure on $E(1)_K$ and we don't know whether $T_C$ is symplectic
with respect to a symplectic structure on $E(1)_K$.
In fact, it is known that $E(1)_K$ doesn't admit any symplectic
structure if the Alexander polynomial of $K$\/ is not monic
\cite{fs:knots}.

On a more positive note, we can easily extend Theorem~\ref{theorem:main}
to $E(n)_K$ $(n\geq 2)$
and more generally to $X_K$, where $X$\/ is a symplectic
4-manifold satisfying
certain topological conditions as in \cite{ep1}.

\begin{proposition}\label{prop:generalization}
Assume that\/ $F$ is a symplectic\/
$2$-torus in a symplectic\/ $4$-manifold\/ $X$.
Suppose that\/ $[F]\in H_2(X)$ is
primitive, $[F]\cdot[F]=0$, and $H^1(X\setminus\nu
F\hspace{1pt})=0$.  If $\,b_2^+(X)=1$, then we also assume
that $\,\overline{SW}
_{\! X\setminus \nu F}\neq 0\,$ and is a finite
sum.
Then there exists
an infinite family of pairwise non-isotopic symplectic tori in\/ $X_K$
representing the homology class\/ $[F]\in H_2(X_K)$ for
any nontrivial fibred knot\/ $K\subset S^3$.
\end{proposition}

%\begin{proof}
The divisibility argument
in the proof of Corollary~\ref{cor:divisibilities}
may not work in this general setting, but after observing
that an isotopy between these tori should preserve $\xi$
and $\tau$, one can resort to a homology basis argument
due to Fintushel and Stern which was announced in \cite{fs:ipam}.
%\end{proof}

It may be possible, as in the rational elliptic surface case,
to show that these non-isotopic tori are inequivalent
under self-diffeomorphisms of $X_K$ once we know the
Seiberg-Witten invariant of\/ $[X \setminus \nu F]$\/ explicitly,
but a general argument
doesn't seem to exist at this moment.

\subsection*{Acknowledgments}
We would like to thank Ronald Fintushel, Sa\v{s}o Strle
and
Stefano Vidussi for their encouragement
and
helpful comments on this and other works of ours.

\end{document}